\documentclass[12pt]{article}

\usepackage{amsmath,amssymb,amsthm,amscd,a4wide,upref}

\usepackage[enableskew]{youngtab}

\usepackage{hyperref}
\hypersetup{colorlinks,citecolor=blue,filecolor=black,linkcolor=blue,urlcolor=blue}
\usepackage{enumerate}
\usepackage{mathtools}
\usepackage{enumitem}

\usepackage{lmodern}     
\usepackage[T1]{fontenc}

\makeatletter
\DeclareTextCommand{\textprime}{\encodingdefault}{%
  \mbox{$\m@th'\kern-\scriptspace$}%
}
\makeatother

\begin{document}


\newcommand{\ad}{{\rm ad}}
\newcommand{\cri}{{\rm cri}}
\newcommand{\row}{{\rm row}}
\newcommand{\col}{{\rm col}}
\newcommand{\End}{{\rm{End}\ts}}
\newcommand{\Rep}{{\rm{Rep}\ts}}
\newcommand{\Hom}{{\rm{Hom}}}
\newcommand{\Mat}{{\rm{Mat}}}
\newcommand{\ch}{{\rm{ch}\ts}}
\newcommand{\chara}{{\rm{char}\ts}}
\newcommand{\diag}{{\rm diag}}
\newcommand{\st}{{\rm st}}
\newcommand{\non}{\nonumber}
\newcommand{\wt}{\widetilde}
\newcommand{\wh}{\widehat}
\newcommand{\ol}{\overline}
\newcommand{\ot}{\otimes}
\newcommand{\la}{\lambda}
\newcommand{\La}{\Lambda}
\newcommand{\De}{\Delta}
\newcommand{\al}{\alpha}
\newcommand{\be}{\beta}
\newcommand{\ga}{\gamma}
\newcommand{\Ga}{\Gamma}
\newcommand{\ep}{\epsilon}
\newcommand{\ka}{\kappa}
\newcommand{\vk}{\varkappa}
\newcommand{\vt}{\vartheta}
\newcommand{\si}{\sigma}
\newcommand{\vs}{\varsigma}
\newcommand{\vp}{\varphi}
\newcommand{\de}{\delta}
\newcommand{\ze}{\zeta}
\newcommand{\om}{\omega}
\newcommand{\Om}{\Omega}
\newcommand{\ee}{\epsilon^{}}
\newcommand{\su}{s^{}}
\newcommand{\hra}{\hookrightarrow}
\newcommand{\ve}{\varepsilon}
\newcommand{\ts}{\,}
\newcommand{\vac}{\mathbf{1}}
\newcommand{\vacu}{|0\rangle}
\newcommand{\di}{\partial}
\newcommand{\qin}{q^{-1}}
\newcommand{\tss}{\hspace{1pt}}
\newcommand{\Sr}{ {\rm S}}
\newcommand{\U}{ {\rm U}}
\newcommand{\BL}{ {\overline L}}
\newcommand{\BE}{ {\overline E}}
\newcommand{\BP}{ {\overline P}}
\newcommand{\AAb}{\mathbb{A}\tss}
\newcommand{\CC}{\mathbb{C}\tss}
\newcommand{\KK}{\mathbb{K}\tss}
\newcommand{\QQ}{\mathbb{Q}\tss}
\newcommand{\SSb}{\mathbb{S}\tss}
\newcommand{\TT}{\mathbb{T}\tss}
\newcommand{\ZZ}{\mathbb{Z}\tss}
\newcommand{\DY}{ {\rm DY}}
\newcommand{\X}{ {\rm X}}
\newcommand{\Y}{ {\rm Y}}
\newcommand{\Z}{{\rm Z}}
\newcommand{\Ac}{\mathcal{A}}
\newcommand{\Lc}{\mathcal{L}}
\newcommand{\Mc}{\mathcal{M}}
\newcommand{\Pc}{\mathcal{P}}
\newcommand{\Qc}{\mathcal{Q}}
\newcommand{\Rc}{\mathcal{R}}
\newcommand{\Sc}{\mathcal{S}}
\newcommand{\Tc}{\mathcal{T}}
\newcommand{\Bc}{\mathcal{B}}
\newcommand{\Ec}{\mathcal{E}}
\newcommand{\Fc}{\mathcal{F}}
\newcommand{\Gc}{\mathcal{G}}
\newcommand{\Hc}{\mathcal{H}}
\newcommand{\Uc}{\mathcal{U}}
\newcommand{\Vc}{\mathcal{V}}
\newcommand{\Wc}{\mathcal{W}}
\newcommand{\Xc}{\mathcal{X}}
\newcommand{\Yc}{\mathcal{Y}}
\newcommand{\Cl}{\mathcal{C}l}
\newcommand{\Ar}{{\rm A}}
\newcommand{\Br}{{\rm B}}
\newcommand{\Ir}{{\rm I}}
\newcommand{\Fr}{{\rm F}}
\newcommand{\Jr}{{\rm J}}
\newcommand{\Or}{{\rm O}}
\newcommand{\GL}{{\rm GL}}
\newcommand{\Spr}{{\rm Sp}}
\newcommand{\Rr}{{\rm R}}
\newcommand{\Zr}{{\rm Z}}
\newcommand{\gl}{\mathfrak{gl}}
\newcommand{\middd}{{\rm mid}}
\newcommand{\ev}{{\rm ev}}
\newcommand{\Pf}{{\rm Pf}}
\newcommand{\Norm}{{\rm Norm\tss}}
\newcommand{\oa}{\mathfrak{o}}
\newcommand{\spa}{\mathfrak{sp}}
\newcommand{\osp}{\mathfrak{osp}}
\newcommand{\f}{\mathfrak{f}}
\newcommand{\g}{\mathfrak{g}}
\newcommand{\h}{\mathfrak h}
\newcommand{\n}{\mathfrak n}
\newcommand{\m}{\mathfrak m}
\newcommand{\z}{\mathfrak{z}}
\newcommand{\Zgot}{\mathfrak{Z}}
\newcommand{\p}{\mathfrak{p}}
\newcommand{\sll}{\mathfrak{sl}}
\newcommand{\gll}{\g^{}_{\ell}}
\newcommand{\gllh}{\wh{\g}^{}_{\ell}}
\newcommand{\gllm}{\g^{}_{\ell,\ell'}}
\newcommand{\glls}{\g^*_{\ell}}
\newcommand{\agot}{\mathfrak{a}}
\newcommand{\bgot}{\mathfrak{b}}
\newcommand{\qdet}{ {\rm qdet}\ts}
\newcommand{\Ber}{ {\rm Ber}\ts}
\newcommand{\HC}{ {\mathcal HC}}
\newcommand{\cdet}{{\rm cdet}}
\newcommand{\rdet}{{\rm rdet}}
\newcommand{\tr}{ {\rm tr}}
\newcommand{\gr}{ {\rm gr}\ts}
\newcommand{\str}{ {\rm str}}
\newcommand{\loc}{{\rm loc}}
\newcommand{\Fun}{{\rm{Fun}\ts}}
\newcommand{\Gr}{{\rm G}}
\newcommand{\sgn}{ {\rm sgn}\ts}
\newcommand{\sign}{{\rm sgn}}
\newcommand{\ba}{\bar{a}}
\newcommand{\bb}{\bar{b}}
\newcommand{\bi}{\bar{\imath}}
\newcommand{\bj}{\bar{\jmath}}
\newcommand{\bk}{\bar{k}}
\newcommand{\bl}{\bar{l}}
\newcommand{\hb}{\mathbf{h}}
\newcommand{\Sym}{\mathfrak S}
\newcommand{\fand}{\quad\text{and}\quad}
\newcommand{\Fand}{\qquad\text{and}\qquad}
\newcommand{\For}{\qquad\text{or}\qquad}
\newcommand{\for}{\quad\text{or}\quad}
\newcommand{\grpr}{{\rm gr}^{\tss\prime}\ts}
\newcommand{\degpr}{{\rm deg}^{\tss\prime}\tss}
\newcommand{\bideg}{{\rm bideg}\ts}

\renewcommand{\theequation}{\arabic{section}.\arabic{equation}}

\numberwithin{equation}{section}

\newtheorem{thm}{Theorem}[section]
\newtheorem{lem}[thm]{Lemma}
\newtheorem{prop}[thm]{Proposition}
\newtheorem{cor}[thm]{Corollary}
\newtheorem{conj}[thm]{Conjecture}
\newtheorem*{mthm}{Main Theorem}
\newtheorem*{mthma}{Theorem A}
\newtheorem*{mthmb}{Theorem B}
\newtheorem*{mthmc}{Theorem C}
\newtheorem*{mthmd}{Theorem D}

\theoremstyle{definition}
\newtheorem{defin}[thm]{Definition}

\theoremstyle{remark}
\newtheorem{remark}[thm]{Remark}
\newtheorem{example}[thm]{Example}
\newtheorem{examples}[thm]{Examples}

\newcommand{\bth}{\begin{thm}}
\renewcommand{\eth}{\end{thm}}
\newcommand{\bpr}{\begin{prop}}
\newcommand{\epr}{\end{prop}}
\newcommand{\ble}{\begin{lem}}
\newcommand{\ele}{\end{lem}}
\newcommand{\bco}{\begin{cor}}
\newcommand{\eco}{\end{cor}}
\newcommand{\bde}{\begin{defin}}
\newcommand{\ede}{\end{defin}}
\newcommand{\bex}{\begin{example}}
\newcommand{\eex}{\end{example}}
\newcommand{\bes}{\begin{examples}}
\newcommand{\ees}{\end{examples}}
\newcommand{\bre}{\begin{remark}}
\newcommand{\ere}{\end{remark}}
\newcommand{\bcj}{\begin{conj}}
\newcommand{\ecj}{\end{conj}}

\newcommand{\bal}{\begin{aligned}}
\newcommand{\eal}{\end{aligned}}
\newcommand{\beq}{\begin{equation}}
\newcommand{\eeq}{\end{equation}}
\newcommand{\ben}{\begin{equation*}}
\newcommand{\een}{\end{equation*}}

\newcommand{\bpf}{\begin{proof}}
\newcommand{\epf}{\end{proof}}

\def\beql#1{\begin{equation}\label{#1}}

\newcommand{\Res}{\mathop{\mathrm{Res}}}

\title{\Large\bf Casimir elements and Sugawara operators for Takiff algebras}

\author{A. I. Molev\footnote{alexander.molev@sydney.edu.au}\\[0.4em]
School of Mathematics and Statistics\\
University of Sydney,
NSW 2006, Australia}

\date{} 
\maketitle


\begin{abstract}
For every simple Lie algebra $\g$ we consider the associated
Takiff algebra $\gll$ defined as the truncated polynomial current Lie algebra
with coefficients in $\g$. We use a matrix presentation of $\gll$ to
give a uniform construction of algebraically independent
generators of the center of the universal enveloping algebra $\U(\gll)$.
A similar matrix presentation for the affine Kac--Moody algebra $\gllh$ is then used
to prove an analogue of the Feigin--Frenkel theorem describing the center
of the corresponding affine vertex algebra at the critical level.
The proof relies on an explicit
construction of a complete set of Segal--Sugawara
vectors for the Lie algebra $\gll$.
\end{abstract}



%

\section{Introduction}
\label{sec:int}

For each simple finite-dimensional
Lie algebra $\g$ over $\CC$ and any positive integer $\ell$ consider
the truncated polynomial current Lie algebra $\gll$ which is defined as the quotient of
$\g\ot\CC[v]$ by the ideal $\g\ot\CC[v]\tss v^{\ell+1}$. The Lie algebra $\gll$
is also called the {\em generalized} or $\ell$-{\em th} {\em Takiff algebra} following the
pioneering work \cite{t:ri}, where such algebras were studied in the case $\ell=1$.
As shown in that paper, the subalgebra of $\gll$-invariants in the symmetric algebra
$\Sr(\gll)$ is an algebra of polynomials. This result was extended by Ra\"{\i}s and
Tauvel~\cite{rt:ip} to all values of $\ell$. More recently, Macedo and Savage~\cite{ms:ip}
proved its multi-parameter generalization, while Panyushev and Yakimova~\cite{py:ta}
showed that this generalization remains valid for a wide class
of Lie algebras $\g$ beyond simple Lie algebras.

The results of Ra\"{\i}s and Tauvel were used by Geoffriau~\cite{g:hh}
to describe properties of an analogue of the Harish-Chandra homomorphism for the
center of the universal enveloping algebra $\U(\gll)$.
Explicit
generators of the center
in type $A$ were given in \cite{m:ce}.
In this case the $\ell$-th Takiff algebra associated with $\gl_n$ coincides with the centralizer
of a certain nilpotent element $e$ in $\gl_{n(\ell+1)}$; namely, $e$ is the direct sum
of $n$ Jordan blocks of size $\ell+1$. The construction of central
elements was extended by Brown and Brundan~\cite{bb:ei}
to arbitrary nilpotents.

Here we give a uniform explicit construction of algebraically independent
generators of the center of $\U(\gll)$ for all simple Lie algebras $\g$ and all $\ell\geqslant 1$.

Then we equip $\gll$ with an invariant symmetric bilinear form by extending
a standard normalized Killing form on $\g$.
The corresponding
affine Kac--Moody algebra $\gllh$ is defined as a central extension of the Lie
algebra of Laurent polynomials $\gll\tss[t,t^{-1}]$.
The vacuum module over $\gllh$ is a vertex algebra whose
center $\z(\gllh)$ is a commutative associative algebra. In the case $\ell=0$
the structure of the center $\z(\wh\g)$
at the critical level
was described by a celebrated theorem of Feigin and Frenkel~\cite{ff:ak} (see also~\cite{f:lc}),
which states that $\z(\wh\g)$ is an algebra of polynomials in infinitely many
variables. We show that this property is shared by the center $\z(\gllh)$
for all $\ell\geqslant 1$.

Our arguments rely on matrix presentations of the Lie algebras $\gll$ and $\gllh$.
Such presentations of the classical Lie algebras and the exceptional Lie algebra of type $G_2$
played a key role in the constructions of generators of the Feigin--Frenkel center
$\z(\wh\g)$ in \cite{m:so} and \cite{mrr:ss}; see also
\cite{o:sf} for a different approach. A recent work by Wendlandt~\cite{w:rm}
provides a significant extension of the matrix techniques
by giving a presentation of $\U(\g)$ for
any simple Lie algebra $\g$ associated with its arbitrary faithful representation.
We recall some results from that paper below
as they will be needed
for our calculations.
To show that our central elements are free generators
of $\z(\gllh)$ we use the classical limit and employ
the Macedo--Savage theorem~\cite{ms:ip} in the particular case of `double' Takiff
algebras.

An analogue of the Feigin--Frenkel theorem for Takiff algebras
in type $A$ was already proved
by Arakawa and Premet~\cite{ap:qm} as a particular case of a more general theorem
describing the centers at the critical level of the affine vertex algebras
associated with centralizers of nilpotent elements in simple Lie algebras.
Explicit generators of the center in type $A$ were produced in \cite{m:cc}.

As in \cite{ap:qm} and \cite{m:cc}, our generators of the center
$\z(\gllh)$ can be used to produce generators of
quantum shift of argument subalgebras of $\U(\gll)$.
We expect that under certain regularity conditions they will be
`quantizations' of the Mishchenko--Fomenko subalgebras of $\Sr(\gll)$
thus yielding a solution of Vinberg's quantization problem
for Takiff algebras.

\section{Matrix presentations}
\label{sec:mp}

We start by recalling some standard tensor product notation.
Any $N\times N$ matrix $X=[X_{ij}]$ with entries in an associative algebra $\Ac$
will be regarded as the element
\beql{mata}
X=\sum_{i,j=1}^N X_{ij}\ot e_{ij}\in\Ac\ot \End\CC^N,
\eeq
where the $e_{ij}\in \End\CC^N$ denote the standard matrix units.
We will need tensor product algebras of the form $\Ac\ot\End(\CC^{N})^{\ot m}$.
For any $a\in\{1,\dots,m\}$ we will denote by $X_a$ the element \eqref{mata}
associated with the $a$-th copy of $\End\CC^{N}$ so that
\ben
X_a=\sum_{i,j=1}^N X_{ij}\ot 1^{\ot(a-1)}\ot e_{ij}\ot 1^{\ot(m-a)}
\in \Ac\ot\End(\CC^{N})^{\ot m}.
\een
Given any element
\ben
C=\sum_{i,j,k,l=1}^N c^{}_{ijkl}\ts e_{ij}\ot e_{kl}\in
\End \CC^N\ot\End \CC^N,
\een
for any two indices $a,b\in\{1,\dots,m\}$ such that $a<b$,
we set
\ben
C_{ab}=\sum_{i,j,k,l=1}^N c^{}_{ijkl}\ts
1^{\ot(a-1)}\ot e_{ij}\ot 1^{\ot(b-a-1)}\ot e_{kl}\ot 1^{\ot(m-b)}\in\End(\CC^{N})^{\ot m}.
\een
Sometimes an additional copy of the endomorphism algebra $\End\CC^N$ labelled
by $0$ will be used so that the notation extends accordingly to that case.

For any $a\in\{1,\dots,m\}$ the {\em partial trace} $\tr_a$
will be understood as the linear map
\ben
\tr_a:\End(\CC^{N})^{\ot m}\to \End(\CC^{N})^{\ot (m-1)}
\een
which acts as the usual trace map on the $a$-th copy of $\End\CC^N$ and
is the identity map on all the remaining copies.
Similarly, the {\em partial transposition} $t_a$ is the linear map on $\End(\CC^{N})^{\ot m}$
which acts as the usual transposition $t:e_{ij}\mapsto e_{ji}$
on the $a$-th copy of $\End\CC^N$ and
is the identity map on all the remaining copies.

For a given simple Lie algebra $\g$ over $\CC$ introduce
a symmetric invariant bilinear form $\langle\ts\ts,\ts\rangle$
as the normalized Killing form
\beql{killi}
\langle X,Y\rangle=\frac{1}{2\tss h^{\vee}}\ts\tr\ts\big(\ad\tss X\ts\ad\tss Y\big),
\eeq
where $h^{\vee}$ is the dual Coxeter number for $\g$.
Fix a basis $J^1,\dots,J^d$ of $\g$ and let
$J_1,\dots,J_d$ be the basis dual to $J^1,\dots,J^d$
with respect to the form \eqref{killi}.
Define Casimir elements by
\beql{omom}
\ol\Om=\sum_{i=1}^d J_i\ot J^i\in\U(\g)\ot\U(\g)\Fand \om=\sum_{i=1}^d J_i\tss J^i\in \U(\g).
\eeq
With the chosen normalization of the Killing form,
the eigenvalue of $\om$ in the adjoint representation equals $2\tss h^{\vee}$
which coincides with the value of the parameter $c_{\g}$ in \cite{w:rm}.
Now let $\rho:\g\to \gl(V)$ be a faithful representation of $\g$ of dimension $\dim V=N$.
Identify the vector space $V$ with $\CC^N$ by choosing a basis and set
\ben
\Om=(\rho\ot\rho)(\ol\Om)\in\End \CC^N\ot\End \CC^N.
\een
Furthermore, define the matrix
\beql{F}
F=\sum_{i,j=1}^N F_{ij}\ot e_{ij}\in \U(\g)\ot\End\CC^N,
\eeq
by setting
$F=-(1\ot \rho)(\ol\Om)$.

The following presentation of $\U(\g)$ is due to Wendlandt~\cite[Proposition~4.4]{w:rm}.

\bpr\label{prop:upres}
The algebra $\U(\g)$ is generated by the elements $F_{ij}$ with $1\leqslant i,j\leqslant N$
subject only to the relations
\beql{ffom}
F_1\tss F_2-F_2\tss F_1=\Om\ts  F_2-F_2\tss \Om
\eeq
in $\U(\g)\ot\End \CC^N\ot\End \CC^N$, and
\beql{symf}
F^2-\big((F^{\tss t})^2\big)^t=h^{\vee}\tss F.
\eeq
\epr

It is easy to see that $[F_1+F_2,\Om]=0$ so that relation \eqref{ffom}
can be written in the equivalent form
\beql{ffomeq}
F_1\tss F_2-F_2\tss F_1=F_1\tss \Om-\Om\ts  F_1.
\eeq

\bre\label{rem:cltypes}
Let $P\in \End \CC^N\ot\End \CC^N$ be the operator permuting the tensor factors
(see also \eqref{p} below).
Proposition~\ref{prop:upres} will hold if
\eqref{symf} is replaced by a system of linear relations
on the entries of the matrix $F$ given by
\beql{symfli}
\tr^{}_2 (\Om\ts  F_2-F_2\tss \Om)P=h^{\vee}\tss F_1.
\eeq
If $\rho$ is the vector representation for the classical types, then in types $B,C$ and $D$
this is equivalent to the skew-symmetry condition
$F+F^{\tss\prime}=0$ with respect to the bilinear
form defining the orthogonal or symplectic Lie algebra. In type $A$ the relation poses
no extra conditions on the generator matrix.
An explicit form of the symmetry conditions \eqref{symfli} for the
$7$-dimensional representation in type $G_2$ can be found in \cite{mrr:ss}.
\ere

One consequence of relation \eqref{ffom} is the following well-known property
of the powers of the generator matrix; cf. \cite[Proposition~4.2.1]{m:so}.

\bco\label{cor:casim}
All elements $\tr\ts F^m$ with $m\geqslant 0$ belong to the center of $\U(\g)$.
\eco

It is also known that by taking $\rho$ to be the lowest-dimension
representation of $\g$, one can choose algebraically independent generators
of the center of $\U(\g)$ among the Casimir elements $\tr\ts F^m$
(with the exception of type $D$, where a Pfaffian-type element $\Pf\ts F$ has to be added).
The required values of $m$ coincide with the degrees of basic invariants
of the symmetric algebra $\Sr(\g)$ as given in
Table~1.


\begin{table}[h!]
\centering
\begin{tabular}{|l|l|}
\hline \ {\rm Type of} $\g$\ \ &\ {\rm Degrees of generators}\ \ \\
\hline
\hline
${A}_n$ & $\ 2$, $3$, $\ldots\,$, $n+1$\\
\hline
${B}_n$ & $\ 2$, $4$, $\ldots\,$,  $2n$\\
\hline
${C}_n$ & $\ 2$, $4$, $\ldots\,$,  $2n$\\
\hline
${D}_n$ & $\ 2$, $4$, $\ldots\,$,  $2n-2$, $n$\\
\hline
${E}_6$ & $\ 2$, $5$, $6$,  $8$, $9$, $12$\\
\hline
${E}_7$ & $\ 2$, $6$, $8$,  $10$, $12$, $14$, $18$\\
\hline
${E}_8$ & $\ 2$, $8$, $12$,  $14$, $18$, $20$, $24$, $30$\\
\hline
${F}_4$ & $\ 2$, $6$, $8$,  $12$\\
\hline
${G}_2$ & $\ 2$, $6$\\
\hline
\end{tabular}

\bigskip

\caption{Degrees of basic invariants} \label{tab:deg}
\end{table}

\noindent
More precisely, the following holds.

\bco\label{cor:genfin}
Except for type $D_n$, the elements
$\tr\ts F^m$ with $m$ running over the values specified in Table~1
are algebraically independent generators of the center of the algebra $\U(\g)$.

In type $D_n$ the elements
$\tr\ts F^m$ with $m=2,4,\dots,2n-2$
and $\Pf\ts F$
are algebraically
independent generators of the center of the algebra $\U(\g)$.
\eco

Corollary~\ref{cor:genfin}
is a classical result for types $A,B,C,D$, but it appears to be less known
for the exceptional types. In those cases, to prove that the elements $\tr\ts F^m$ are algebraically
independent generators, one only needs to verify that their
top degree components are basic $\g$-invariants of the symmetric algebra $\Sr(\g)$.
The latter property
goes back to Kuin\textprime~\cite{k:pp}, \cite{k:ppc}, with some cases previously considered by
Coxeter~\cite{c:pg} ($E_6$) and
Takeuchi~\cite{t:pc} $(F_4)$.
A direct claim about
these Casimir elements was made in
\cite{b:cc}, \cite{bb:co}.

To give more details for the exceptional types, recall that
the theorem of Kuin\textprime{}, whose proof is given
in \cite[\S2.2]{k:ppc}, reads as follows. Suppose that
$\La_1,\dots,\La_N$ are the weights of the lowest-dimension
representation $(V,\rho)$ of $\g$. The weights are understood as
elements of the Cartan subalgebra $\h$ of $\g$, where
$\h^*$ and $\h$ are identified via the form \eqref{killi}.
That is, $V$ has a basis $v_1,\dots,v_N$ such that any element $X\in\h$ acts by
\ben
X\tss v_a=\langle\La_a,X\rangle\tss v_a,\qquad a=1,\dots,N.
\een
The theorem states that the power sums
\beql{pm}
P_m=\sum_{a=1}^N \La_a^m
\eeq
with $m$ running over the respective degrees in Table 1, are algebraically
independent generators of the subalgebra of $W$-invariants $\Sr(\h)^W$ in $\Sr(\h)$, where
$W$ denotes the Weyl group of $\g$.

Returning to Corollary~\ref{cor:genfin},
choose a special form of the Casimir element $\ol\Om$ in \eqref{omom}
by taking
a basis $J^1,\dots,J^d$ of $\g$ such that
$J^1,\dots,J^n$ form a basis of the Cartan subalgebra $\h$, while
the remaining $J^i$s are root vectors. Then
the vectors $J_1,\dots,J_n$ of the dual basis will also belong to $\h$.
Consider the entries $F_{ij}$ of the matrix $F$ defined in \eqref{F}
as elements of the symmetric algebra $\Sr(\g)$. A nonzero contribution
to the image of the element $\tr\ts F^m\in \Sr(\g)^{\g}$ under the
Chevalley isomorphism $\Sr(\g)^{\g}\to \Sr(\h)^W$ will only come from
its part $\tr\ts \ol F^{\ts m}$, where
\ben
\ol F=-\sum_{i=1}^n J_i\ot \rho(J^i).
\een
By our choice of parameters, $\rho(J^i)$ is a diagonal matrix of the form
\ben
\rho(J^i)=\sum_{a=1}^N \langle\La_a,J^i\rangle\ts e_{aa}.
\een
Hence, the Chevalley image of $\tr\ts F^m$ equals
\ben
\sum_{a=1}^N\Big({-}\sum_{i=1}^n \langle\La_a,J^i\rangle\ts J_i\Big)^m
=(-1)^m\sum_{a=1}^N \La_a^m=(-1)^m P_m
\een
which thus coincides with the element \eqref{pm}, up to a sign.

\section{Casimir elements for Takiff algebras}
\label{sec:cetakiff}

We will use the presentation of $\g$ associated with an arbitrary
faithful representation $\rho$,
as given in
Proposition~\ref{prop:upres}. Introduce elements of the Takiff algebra $\gll$
as defined in the Introduction, by $F_{ij}^{(r)}=F_{ij} v^r$ for $r=0,1,\dots,\ell$.
We combine them into the respective matrices
\ben
F^{(r)}=\sum_{i,j=1}^N F_{ij}^{(r)}\ot e_{ij}\in \U(\gll)\ot \End\CC^N.
\een
For a variable $u$ set
\ben
F(u)=F^{(0)}+F^{(1)}u+\dots+F^{(\ell)}u^{\ell}.
\een
The traces of powers of this matrix are polynomials in $u$ of the form
\ben
\tr\ts F(u)^m=\sum_r \theta_m^{\tss(r)}\tss u^r,\qquad \theta_m^{\tss(r)}\in \U(\gll).
\een

\bpr\label{prop:casita}
All elements $\theta_m^{\tss(r)}$
with $m\geqslant 1$ and $r=m\ell,m\ell-1,\dots,m\ell-\ell$
belong to the center of the algebra $\U(\gll)$.
\epr

\bpf
Since $\g=[\g,\g]$, the Lie algebra $\gll$ is generated by the elements $F_{ij}^{(0)}$
and $F_{ij}^{(1)}$ with $1\leqslant i,j\leqslant N$.
Hence, it is sufficient
to verify that all commutators
\ben
\big[F_{ij}^{(0)}, \tr\ts F(u)^m\big]\Fand \big[F_{ij}^{(1)}, \tr\ts F(u)^m\big]
\een
are polynomials in $u$ whose degrees are less than $m\ell-\ell$.
We will use the matrix
notation of Section~\ref{sec:mp} and consider the algebra
$\U(\gll)\ot\End \CC^N\ot\End \CC^N$ with the tensor factors $\End \CC^N$ labelled by $0$ and $1$.
Relation \eqref{ffom} implies
\begin{multline}
\big[F^{(0)}_0, F(u)_1^m\big]=\sum_{i=1}^m F(u)_1^{\tss i-1}
\big[F^{(0)}_0,F(u)_1\big]\tss F(u)_1^{\tss m-i}\\
=\sum_{i=1}^m F(u)_1^{\tss i-1} \big[\Om_{\tss 01},F(u)_1\big]\tss F(u)_1^{\tss m-i}
=\big[\Om_{\tss 01},F(u)_1^{\tss m}\big].
\non
\end{multline}
Therefore, taking the partial trace $\tr^{}_1$ on both
sides and using its cyclic property we can conclude that
$[F^{(0)}_0, \tr_1\ts F(u)_1^m]=0$. We also have
\ben
u\ts\big[F^{(1)}_0, F(u)_1\big]=\big[\Om_{\tss 01},F(u)_1-F^{(0)}_1\big].
\een
Hence, a similar calculation gives
\ben
u\ts\big[F^{(1)}_0, \tr_1\ts F(u)_1^m\big]
=-\sum_{i=1}^m \tr_1\ts F(u)_1^{\tss i-1} \big[\Om_{\tss 01},F^{(0)}_1\big]\tss F(u)_1^{\tss m-i}.
\een
However, the degree of the polynomial in $u$ on the right hand side is $m\ell-\ell$
so that the commutator $[F^{(1)}_0, \tr_1\ts F(u)_1^m]$
is a polynomial of degree less than $m\ell-\ell$, as required.
\epf

\bre\label{rem:gln}
Proposition~\ref{prop:casita} and its proof extend
to the reductive Lie algebra $\gl_N$, where the $F_{ij}$ should be replaced
by the standard basis elements $E_{ij}$ of $\gl_N$.
If $\rho$ is the vector representation
of $\gl_N$, the Casimir elements provided by the proposition
are closely related to one of the families produced in
\cite{m:ce}.
\qed
\ere

Now suppose that $\g=\oa_{2n}$ is the orthogonal Lie algebra (of type $D_n$).
We will use its presentation where the elements of $\g$ are skew-symmetric $2n\times 2n$ matrices
with the usual matrix commutator. In relation \eqref{ffom} we have
\ben
\Om=\sum_{i,j=1}^{2n} (e_{ij}\ot e_{ji}-e_{ij}\ot e_{ij}),
\een
while \eqref{symfli} is equivalent to $F+F^t=0$.
Define the {\em Pfaffian} of the matrix $F(u)$ by the formula
\beql{pfcanon}
\Pf\ts F(u)=\sum_{\si}\sgn\si\cdot
F_{\si(1)\ts\si(2)}(u)\dots F_{\si(2n-1)\ts\si(2n)}(u),
\eeq
summed over the elements $\si$ of the subset $\Ac_{2n}\subset \Sym_{2n}$
of the symmetric group $\Sym_{2n}$
which consists of the permutations with the properties
$\si(2k-1)<\si(2k)$ for all $k=1,\dots,n$ and $\si(1)<\si(3)<\dots<\si(2n-1)$.
Introduce the coefficients of the polynomial $\Pf\ts F(u)$ by
\ben
\Pf\ts F(u)=\sum_r \pi^{(r)}u^r.
\een

\bpr\label{prop:pfa}
All coefficients $\pi^{(r)}$ with $r=n\ell,n\ell-1,\dots,n\ell-\ell$
belong to the center of the algebra $\U(\gll)$.
\epr

\bpf
It is enough to prove that the commutators
\beql{annihpf}
\big[F_{ij}^{(0)}, \Pf\ts F(u)\big]
\Fand \big[F_{ij}^{(1)}, \Pf\ts F(u)\big]
\eeq
are polynomials in $u$ of degree less than $n\ell-\ell$.
Note that for any permutation $\pi\in\Sym_{2n}$ the mapping
$
F_{ij}^{(r)}\mapsto F_{\pi(i)\ts\pi(j)}^{(r)}
$
defines an automorphism of the Lie algebra $\gll$.
Hence, it is sufficient to verify the required
properties of the commutators in \eqref{annihpf}
for $i=1$ and $j=2$.
This follows by a straightforward calculation with the use
of the commutation relations
\ben
\big[F_{ij}^{(0)},F_{kl}(u)\big]
=\de_{kj}\ts
F_{il}(u)-\de_{il}\ts F_{kj}(u)
-\de_{ki}\ts F_{jl}(u)+\de_{jl}\ts F_{ki}(u)
\een
and
\ben
u\big[F_{ij}^{(1)},F_{kl}(u)\big]
=\big[F_{ij}^{(0)},F_{kl}(u)-F_{kl}^{(0)}\big]
\een
and arguing as in the proof of Proposition~\ref{prop:casita}.
\epf

\bre\label{rem:pfa}
Analogous expressions for the Pfaffian-type central elements can
also be written for the presentations
of the orthogonal Lie algebra associated with arbitrary non-degenerate
symmetric bilinear forms on $\CC^{2n}$; cf. \cite[Sec.~8]{m:so}.
\qed
\ere

Return to an arbitrary simple Lie algebra $\g$ and suppose now that $\rho$
is the lowest-dimension representation of $\g$.

\bth\label{thm:gencas}
Except for type $D_n$, the elements
$\theta_m^{\tss(r)}$ with $m$ running over the values specified in Table~1
and $r=m\ell,m\ell-1,\dots,m\ell-\ell$, are algebraically
independent generators of the center of the algebra $\U(\gll)$.

In type $D_n$ the elements
$\theta_m^{\tss(r)}$ with $m=2,4,\dots,2n-2$
and $r=m\ell,m\ell-1,\dots,m\ell-\ell$ together with
$\pi^{(r)}$ with $r=n\ell,n\ell-1,\dots,n\ell-\ell$,
are algebraically
independent generators of the center of the algebra $\U(\gll)$.
\eth

\bpf
All these elements belong to the center of $\U(\gll)$ by Propositions~\ref{prop:casita}
and \ref{prop:pfa}. Their symbols in the symmetric algebra $\Sr(\gll)$ are
$\gll$-invariants. Moreover, these invariants are associated with basic
$\g$-invariants in $\Sr(\g)$ in the way described in \cite[Sec.~3.1]{rt:ip}.
Therefore, applying \cite[Th\'{e}or\`{e}me~4.5]{rt:ip}, we can conclude that
the symbols are algebraically independent generators of the subalgebra
of $\gll$-invariants in $\Sr(\gll)$. This implies the desired property
of the central elements in $\U(\gll)$.
\epf

\section{Segal--Sugawara vectors}
\label{sec:ssvec}

We will identify the Lie algebra $\g$ with a subalgebra of $\gll$
via the embedding $F_{ij}\mapsto F^{(0)}_{ij}$. Extend the form \eqref{killi}
defined on this subalgebra to the Lie algebra $\gll$ by positing that
all elements $F^{(r)}_{ij}$ with $r=1,\dots,\ell$ belong to its kernel.
This defines a symmetric invariant bilinear form $\langle\ts\ts,\ts\rangle$
on $\gll$. The corresponding affine Kac--Moody algebra $\gllh$
is the central
extension
\ben
\gllh=\gll\tss[t,t^{-1}]\oplus\CC K,
\een
where $\gll[t,t^{-1}]$ is the Lie algebra of Laurent
polynomials in $t$ with coefficients in $\gll$. For any $r\in\ZZ$ and $X\in\gll$
we will write $X[r]=X\ts t^r$. The commutation relations of the Lie algebra $\gllh$
have the form
\beql{affrel}
\big[X[r],Y[s]\big]=[X,Y][r+s]+r\ts\de_{r,-s}\langle X,Y\rangle\ts K,
\qquad X, Y\in\gll,
\eeq
and the element $K$ is central in $\gllh$.

The {\em vacuum module at the level} $k\in\CC$
over $\gllh$
is the quotient
$
V_k(\gll)=\U(\gllh)/\Ir,
$
where $\Ir$ is the left ideal of $\U(\gllh)$ generated by $\gll[t]$
and the element $K-k$.
The Poincar\'e--Birkhoff--Witt theorem implies that this quotient
is isomorphic to the universal enveloping algebra
$\U\big(t^{-1}\gll[t^{-1}]\big)$, as a vector space.
The vacuum module is equipped with a vertex algebra structure; see e.g.
\cite{f:lc}, \cite{k:va}.
We will call the level $k=-(\ell+1)\tss h^{\vee}$ {\em critical},
as the vacuum module $V_{\cri}(\gll)$ at this level turns out to exhibit similar properties
to its counterpart for $\ell=0$.
We will denote by $\z(\gllh)$ the {\em center} of the vertex algebra $V_{\cri}(\gll)$
which is defined as the subspace
\ben
\z(\gllh)=\{v\in V_{\cri}(\gll)\ |\ \gll[t]\tss v=0\}.
\een
It follows from the axioms of vertex algebra that $\z(\gllh)$
is a unital commutative associative algebra which can be regarded as a
subalgebra of $\U\big(t^{-1}\gll[t^{-1}]\big)$.
This subalgebra is invariant with respect to the
{\em translation operator}
$T$ which is
the derivation of the algebra $\U\big(t^{-1}\gll[t^{-1}]\big)$
whose action on the generators is given by
\beql{tran}
T:X[r]\mapsto -r\tss X[r-1],\qquad X\in\gll, \quad r<0.
\eeq

Any element of $\z(\gllh)$
is called a {\em Segal--Sugawara vector\/}.
By the Feigin--Frenkel theorem \cite{ff:ak}, \cite{f:lc},
in the case $\ell=0$ the center $\z(\wh\g)$ contains
a {\em complete set of Segal--Sugawara vectors} $S_1,\dots,S_n$,
which means that the translations $T^r S_p$ with $r\geqslant 0$ and $p=1,\dots,n({={\rm rank}}\ts\g)$
are algebraically independent generators of the algebra $\z(\wh\g)$.

Our goal is to prove that this property is shared by the algebras $\z(\gllh)$
for all $\ell\geqslant 1$. Note that this has already been proved in type $A$
by Arakawa and Premet~\cite{ap:qm} as a particular case of a more general theorem
on affine vertex algebras associated with centralizers of nilpotent elements in $\g$;
see also \cite{m:cc} for an explicit construction of
a complete set of Segal--Sugawara vectors.

We begin by producing some families of Segal--Sugawara vectors for $\gll$
and then will show how to choose a complete set of such vectors.
As in Section~\ref{sec:cetakiff},
we will use the presentation of $\g$ associated with an arbitrary
faithful representation $\rho$, given in
Proposition~\ref{prop:upres}. Introduce polynomials in $u$ of the form
\ben
\Fc(u)=F^{(0)}[-1]+F^{(1)}[-1]\tss u+\dots+F^{(\ell)}[-1]\tss u^{\ell},
\een
where $F^{(r)}[p]$ denotes the matrix
\beql{matr}
F^{(r)}[p]=\sum_{i,j=1}^N F_{ij}^{(r)}[p]\ot e_{ij},\qquad p\in\ZZ.
\eeq
Define elements $\Theta_m^{(r)}\in V_{\cri}(\gll)\cong \U\big(t^{-1}\gll[t^{-1}]\big)$ as the coefficients
of the polynomial
\ben
\tr\ts \Fc(u)^m=\sum_r \Theta_m^{(r)}\tss u^r.
\een

\bpr\label{prop:ssita}
Suppose that $\ell\geqslant 1$. Then
all coefficients $\Theta_m^{(r)}$
with the parameters $m\geqslant 1$ and $r=m\ell,m\ell-1,\dots,m\ell-\ell$
are Segal--Sugawara vectors for $\gll$.
\epr

\bpf
Since $F_{ij}^{(0)}[0], F_{ij}^{(0)}[1]$
and $F_{ij}^{(1)}[0]$ with $1\leqslant i,j\leqslant N$
are generators of the Lie algebra $\gll[t]$,
it is sufficient
to verify that they annihilate the elements $\Theta_m^{(r)}$.
Using the matrix notation as in the proof
of Proposition~\ref{prop:casita}, we derive from
\eqref{ffom} that
\ben
\big[F^{(0)}[0]_0, \Fc(u)_1\big]=\big[\Om_{\tss 01},\Fc(u)_1\big].
\een
Hence, $F^{(0)}[0]_0\ts\tr_1\ts \Fc(u)_1^m=0$ in $V_{\cri}(\gll)$
which follows in the same way as
in the proof
of Proposition~\ref{prop:casita}.
As in that proof, we also have
\ben
u\ts\big[F^{(1)}[0]_0, \Fc(u)_1\big]=\big[\Om_{\tss 01},\Fc(u)_1-F^{(0)}[-1]_1\big],
\een
which implies that the degree of the polynomial $F^{(1)}[0]_0\ts\tr_1\ts \Fc(u)_1^m$ in $u$
is less than $m\ell-\ell$.

Furthermore, as an immediate consequence
of \eqref{ffomeq} and \eqref{affrel}, we obtain
\ben
\big[F^{(0)}[r]_0,F^{(0)}[s]_1\big]=\big[F^{(0)}[r+s]_0,\Om_{\tss 01}\big]+
r\ts\de_{r,-s}\tss \Om_{\tss 01}\tss K.
\een
Hence, for the remaining generators of $\gll[t]$ we have
\ben
\big[F^{(0)}[1]_0, \Fc(u)_1\big]=\big[\Fc(0,u)_0,\Om_{\tss 01}\big]+\Om_{\tss 01}\tss K,
\een
where we used the notation
\beql{fou}
\Fc(0,u)=F^{(0)}[0]+F^{(1)}[0]\tss u+\dots+F^{(\ell)}[0]\tss u^{\ell}.
\eeq
Calculating in the vacuum module
we then find
\ben
F^{(0)}[1]_0\ts\tr^{}_{1}\ts \Fc(u)_1^m=
\sum_{i=1}^m \tr^{}_{1}\ts \Fc(u)_1^{i-1}\big({-}\Om_{\tss 01}\tss \Fc(0,u)_0+\Fc(0,u)_0\tss
\Om_{\tss 01}+\Om_{\tss 01}\tss K\big)
\Fc(u)_1^{m-i}.
\een
Observe that modulo a polynomial in $u$ of degree less than $\ell$, we can write
\beql{cofoure}
\big[\Fc(0,u)_0,\Fc(u)_1\big]\equiv
(\ell+1)\tss\big[\Om_{\tss 01}, F^{(\ell)}[-1]_1\tss u^{\ell}\big]
\equiv (\ell+1)\tss\big[\Om_{\tss 01}, \Fc(u)_1\big].
\eeq
Therefore,
\beql{fofki}
\Fc(0,u)_0\tss \Fc(u)_1^{m-i}\equiv
(\ell+1)\tss\big[\Om_{\tss 01},\tss \Fc(u)_1^{m-i}\big]
\eeq
modulo a polynomial in $u$ of degree less than $(m-i)\tss\ell$.

Now use the following general property of the partial transposition $t_1$:
\beql{sklyanin-lemma}
\tr^{}_{1}\ts X\tss Y=\tr^{}_{1}\ts X^{t_1}\tss Y^{t_1}.
\eeq
Taking $X=\Fc(u)_1^{i-1}\tss \Fc(0,u)_0$ and $Y=\Om_{\tss 01}\tss \Fc(u)_1^{m-i}$ we obtain
\ben
\tr^{}_{1}\ts \Fc(u)_1^{i-1}\tss \Fc(0,u)_0\tss\Om_{\tss 01}\tss \Fc(u)_1^{m-i}=
\tr^{}_{1}\ts \big(\Fc(u)^{i-1}\big)^t_1\tss \Fc(0,u)_0\tss
\big(\Fc(u)^{m-i}\big)^t_1\tss(\Om_{\tss 01})^{t_1}.
\een
Applying the transposition to both sides of \eqref{fofki} we get
\ben
\Fc(0,u)_0\tss
\big(\Fc(u)^{m-i}\big)^t_1\equiv (\ell+1)\tss\big[\big(\Fc(u)^{m-i}\big)^t_1,(\Om_{\tss 01})^{t_1}\big].
\een
Thus, bringing the calculations together, we obtain the following
relation in $V_{\cri}(\gll)$ modulo
a polynomial in $u$ of degree less than $m\tss \ell-\ell$:
\begin{multline}
F^{(0)}[1]_0\ts\tr^{}_{1}\ts \Fc(u)_1^m\equiv
K\tss \sum_{i=1}^m\tr^{}_{1}
\ts\Fc(u)_1^{i-1}\tss \Om_{\tss 01}\tss \Fc(u)_1^{m-i}\\
{}+(\ell+1)\sum_{i=1}^m\tr^{}_{1}
\ts\Big(\Fc(u)_1^{i-1}\tss \Om_{\tss 01}\tss \Fc(u)_1^{m-i}\tss \Om_{\tss 01}
-\Fc(u)_1^{i-1}\Om_{\tss 01}^{\tss 2}\tss \Fc(u)_1^{m-i}\\[0.4em]
{}-\big(\Fc(u)^{i-1}\big)^t_1\tss(\Om_{\tss 01})^{t_1}\tss
\big(\Fc(u)^{m-i}\big)^t_1\tss(\Om_{\tss 01})^{t_1}
+\big(\Fc(u)^{i-1}\big)^t_1\tss\big(\Fc(u)^{m-i}\big)^t_1\tss\big((\Om_{\tss 01})^{t_1}\big)^2\Big).
\non
\end{multline}
The application of \eqref{sklyanin-lemma} to the terms in the last line
brings this expression to the form
\begin{multline}
F^{(0)}[1]_0\ts\tr^{}_{1}\ts \Fc(u)_1^m\equiv
\sum_{i=1}^m\tr^{}_{1}
\ts\Fc(u)_1^{i-1}\tss \Big(K\tss\Om_{\tss 01}+(\ell+1)
\Big(\big((\Om_{\tss 01})^{t_1}\big)^2\Big)^{t_1}-(\ell+1)\Om_{\tss 01}^{\tss 2}\Big)
\tss \Fc(u)_1^{m-i}\\[0.3em]
{}+(\ell+1)\sum_{i=1}^m\tr^{}_{1}
\ts\Big(\Fc(u)_1^{i-1}\tss \Om_{\tss 01}\tss \Fc(u)_1^{m-i}\tss \Om_{\tss 01}
-\Om_{\tss 01}\tss\Fc(u)_1^{i-1}\tss \Om_{\tss 01}\tss \Fc(u)_1^{m-i}\Big).
\non
\end{multline}
Recall that $K=-(\ell+1)\tss h^{\vee}$ at the critical level and so
the first sum is zero. This follows from
the identity
\ben
\Om_{\tss 01}^{\tss 2}-\Big(\big((\Om_{\tss 01})^{t_1}\big)^2\Big)^{t_1}+h^{\vee}\tss\Om_{\tss 01}=0
\een
which is a consequence of relation \eqref{symf};
we just need to apply $\rho\ot 1$ to its both sides and note that
$\Om_{\tss 01}=-(\rho\ot 1)(F)$.

The second sum is a polynomial in $u$ of degree at most $m\tss \ell-\ell$.
It remains to verify that the coefficient of $u^{m\tss \ell-\ell}$ in the sum is zero.
This coefficient equals
\beql{coeffi}
\sum_{i=1}^m\tr^{}_{1}
\ts\Big(\Phi_1^{i-1}\tss \Om_{\tss 01}\tss \Phi_1^{m-i}\tss \Om_{\tss 01}
-\Om_{\tss 01}\tss\Phi_1^{i-1}\tss \Om_{\tss 01}\tss \Phi_1^{m-i}\Big),
\eeq
where we set $\Phi=F^{(\ell)}[-1]$ for brevity.
In the algebra
\ben
\U(\gllh)\ot\End \CC^N\ot\End \CC^N\ot \End \CC^N
\een
with the tensor factors $\End \CC^N$ labelled by $0$, $1$ and $2$ we can write
\ben
\tr^{}_{1}\ts\Om_{\tss 01}\tss\Phi_1^{i-1}\tss \Om_{\tss 01}\tss \Phi_1^{m-i}
=\tr^{}_{1,2}\ts \Om_{\tss 01}\tss\Phi_1^{i-1}\tss \Om_{\tss 01}\tss \Phi_2^{m-i}P_{12},
\een
where
\beql{p}
P_{12}=\sum_{i,j=1}^N e_{ij}\ot e_{ji}
\eeq
is the permutation operator. We have used the relation $\Phi_2^{m-i}P_{12}=P_{12}\tss\Phi_1^{m-i}$
and observed that $\tr_2\ts P_{12}=1$. Since $\ell\geqslant 1$, the
matrix elements of the matrix $\Phi$ pairwise commute and so
\ben
\tr^{}_{1,2}\ts \Om_{\tss 01}\tss\Phi_1^{i-1}\tss \Om_{\tss 01}\tss \Phi_2^{m-i}P_{12}
=\tr^{}_{1,2}\ts \Phi_2^{m-i}\tss \Om_{\tss 01}\tss\Phi_1^{i-1}\tss \Om_{\tss 01}\tss P_{12}
=\tr^{}_{1,2}\ts \Phi_2^{m-i}\tss P_{12}\tss \Om_{\tss 02}\tss\Phi_2^{i-1}\tss \Om_{\tss 02}
\een
which equals
\ben
\tr^{}_{2}\ts \Phi_2^{m-i}\tss \Om_{\tss 02}\tss\Phi_2^{i-1}\tss \Om_{\tss 02}
=\tr^{}_{1}\ts \Phi_1^{m-i}\tss \Om_{\tss 01}\tss\Phi_1^{i-1}\tss \Om_{\tss 01}.
\een
Thus, the expression \eqref{coeffi} is zero so that
$F^{(0)}[1]_0\ts\tr^{}_{1}\ts \Fc(u)_1^m$ is a polynomial in $u$ of degree less than
$m\tss \ell-\ell$.
\epf

\bre\label{rem:lze}
Proposition~\ref{prop:ssita} does not hold under the assumption $\ell=0$;
see e.g. \cite[Secs~7.1, 8.2 and 8.4]{m:so} for counterexamples in classical types.
\qed
\ere

Consider now the orthogonal Lie algebra $\g=\oa_{2n}$
with its presentation used in Section~\ref{sec:cetakiff}.
Define elements $\Pi^{(r)}\in V_{\cri}(\gll)$ as the coefficients
of the {\em Pfaffian} of the matrix $\Fc(u)$
\ben
\Pf\ts \Fc(u)=\sum_r \Pi^{(r)}u^r,
\een
where $\Pf\ts \Fc(u)$ is given by formula \eqref{pfcanon}
applied to the matrix $\Fc(u)$.

\bpr\label{prop:pfaaff}
All coefficients $\Pi^{(r)}$ with $r=n\ell,n\ell-1,\dots,n\ell-\ell$
are Segal--Sugawara vectors for $\gll$.
\epr

\bpf
Given any permutation $\pi\in\Sym_{2n}$, the mapping
\ben
F^{(r)}_{ij}[p]\mapsto F^{(r)}_{\pi(i)\ts\pi(j)}[p],\qquad K\mapsto K
\een
defines an automorphism of the Lie algebra $\gllh$.
Therefore, it is sufficient
to verify that the elements
$F_{12}^{(0)}[0], F_{12}^{(1)}[0]$
and $F_{12}^{(0)}[1]$
acting in the vacuum module $V_{\cri}(\gll)$
annihilate the given coefficients $\Pi^{(r)}$.
The argument for the first two elements is
a straightforward calculation with the use
of the commutation relations involving the entries of the matrix $\Fc(u)=[\Fc_{ij}(u)]$,
\ben
\big[F_{ij}^{(0)}[0],\Fc_{kl}(u)\big]
=\de_{kj}\ts
\Fc_{il}(u)-\de_{il}\ts \Fc_{kj}(u)
-\de_{ki}\ts \Fc_{jl}(u)+\de_{jl}\ts \Fc_{ki}(u)
\een
and
\ben
u\big[F_{ij}^{(1)}[0],\Fc_{kl}(u)\big]
=\big[\Fc_{ij}^{(0)},\Fc_{kl}(u)-F_{kl}^{(0)}[-1]\big];
\een
cf. the case $\ell=0$ in \cite[Proposition~8.1.4]{m:so}.
For the remaining element we use the relations
\ben
\bal
\big[F_{ij}^{(0)}[1],\Fc_{kl}(u)\big]
=\de_{kj}\ts
\Fc_{il}(0,u)-\de_{il}\ts \Fc_{kj}(0,u)
-\de_{ki}\ts \Fc_{jl}(0,u)&{}+\de_{jl}\ts \Fc_{ki}(0,u)\\
{}&+K\big(\de_{kj}\ts\de_{il}-\de_{ki}\ts\de_{jl}\big)
\eal
\een
involving the entries of the matrix $\eqref{fou}$.
Consider first the summands in formula \eqref{pfcanon} for the matrix $\Fc(u)$ with
$\si(1)=1$ and $\si(2)=2$. In the vacuum module we have
\ben
F^{(0)}_{12}[1]\ts \Fc_{1\ts 2}(u)\tss \Fc_{\si(3)\ts \si(4)}(u)
\dots \Fc_{\si(2n-1)\ts\si(2n)}(u)=-K\ts \Fc_{\si(3)\ts \si(4)}(u)
\dots \Fc_{\si(2n-1)\ts\si(2n)}(u).
\een
Furthermore, let $\tau\in\Ac_{2n}$
with $\tau(2)>2$. Then $\tau(3)=2$ and $\tau(4)>2$ and so
\ben
F^{(0)}_{12}[1]\ts \Fc_{1\ts \tau(2)}(u)\tss \Fc_{2\ts \tau(4)}(u)
\dots \Fc_{\tau(2n-1)\ts\tau(2n)}(u)
=-\Fc_{2\ts \tau(2)}(0,u)\tss \Fc_{2\ts \tau(4)}(u)\dots \Fc_{\tau(2n-1)\ts\tau(2n)}(u).
\een
Considering this expression modulo a polynomial in $u$ of degree less that
$n\ell-\ell$, we may apply relations
\eqref{cofoure} to conclude that the expression coincides with
\ben
(\ell+1)\ts \Fc_{\tau(2)\ts \tau(4)}(u)\dots \Fc_{\tau(2n-1)\ts\tau(2n)}(u).
\een
Suppose now that $\si\in\Ac_{2n}$ in an element with
the fixed values $\si(1)=1$ and $\si(2)=2$, and calculate
the coefficient of the monomial
\ben
\Fc_{\si(3)\ts \si(4)}(u)\dots \Fc_{\si(2n-1)\ts\si(2n)}(u)
\een
in the expansion of $F^{(0)}_{12}[1]\ts \Pf\ts \Fc(u)$
modulo a polynomial in $u$ of degree less that
$n\ell-\ell$. Essentially the same calculation was already performed
in the case $\ell=0$ in \cite[Proposition~8.1.4]{m:so} showing that
the coefficient equals
${-}K-(\ell+1)(2n-2)$, which is zero at the critical level since $h^{\vee}=2n-2$.
\epf

Now let $\g$ be an arbitrary simple Lie algebra and suppose that $\rho$
is the lowest-dimension representation of $\g$.
The next theorem shows that a natural extension of the Feigin--Frenkel theorem~\cite{ff:ak}
holds for the Lie algebra $\gll$ with $\ell\geqslant 1$; that is, the center $\z(\gllh)$
is an algebra of polynomials.

\bth\label{thm:gencasaff}
Except for type $D_n$, the elements
$\Theta_m^{(r)}$ with $m$ running over the values specified in Table~1
and $r=m\ell,m\ell-1,\dots,m\ell-\ell$, form a complete set of Segal--Sugawara vectors
for the Lie algebra $\gll$.

In type $D_n$ the elements
$\Theta_m^{(r)}$ with $m=2,4,\dots,2n-2$
and $r=m\ell,m\ell-1,\dots,m\ell-\ell$ together with
$\Pi^{(r)}$ with $r=n\ell,n\ell-1,\dots,n\ell-\ell$,
form a complete set of Segal--Sugawara vectors
for the Lie algebra $\gll$.
\eth

\bpf
All elements are Segal--Sugawara vectors by Propositions~\ref{prop:ssita}
and \ref{prop:pfaaff}. We need to show that
all the shifted vectors $T^s\tss\Theta_m^{(r)}$
(together with $T^s\tss\Pi^{(r)}$ in type $D_n$) with $s\geqslant 0$
are algebraically
independent generators of the algebra $\z(\gllh)$. We will follow the approach
which was used for the case
$\ell=0$; cf.~\cite[Secs~3.3 and 3.4]{f:lc} and \cite[Sec.~6.3]{m:so}.

Regard $V_{\cri}(\gll)$ as a $\gll[t]$-module obtained
by restriction of the action of $\gllh$ to the subalgebra $\gll[t]$.
By identifying the vector space $V_{\cri}(\gll)$ with the algebra
$\U\big(t^{-1}\gll[t^{-1}]\big)$ and using its canonical filtration, equip
the associated graded space $\gr V_{\cri}(\gll)$
with the structure of
a $\gll[t]$-module. As a vector space, $\gr V_{\cri}(\gll)$ will be identified
with the symmetric algebra $\Sr\big(t^{-1}\gll[t^{-1}]\big)$.
The action of $\gll[t]$ on $\Sr\big(t^{-1}\gll[t^{-1}]\big)$ is obtained
by extending the adjoint representation of $\gll[t]$ on $\gll[t,t^{-1}]/\gll[t]\cong t^{-1}\gll[t^{-1}]$
to the symmetric algebra.

The symbols $\ol\Theta_m^{\tss(r)}$ of the Segal--Sugawara vectors $\Theta_m^{(r)}$
(and the symbols $\ol\Pi^{\ts(r)}$ of the Segal--Sugawara vectors $\Pi^{(r)}$ in type $D_n$)
belong to the subalgebra
$\Sr\big(t^{-1}\gll[t^{-1}]\big)^{\gll[t]}$
of $\gll[t]$-invariants in $\Sr\big(t^{-1}\gll[t^{-1}]\big)$.
The translation operator defined in \eqref{tran} induces a derivation of the symmetric
algebra which we will also denote by $T$.
We only need to verify that all the shifts
$T^s\tss\ol\Theta_m^{\ts(r)}$
(together with $T^s\tss\ol\Pi^{\ts(r)}$ in type $D_n$) with $s\geqslant 0$
are algebraically
independent generators of the algebra of invariants.

For the rest of the proof we will work with the symmetric algebra and so we will
keep the same notation $F_{ij}^{(r)}[-p-1]$ with $p\geqslant 0$ for the images of these elements
of the enveloping algebra in the associated graded
algebra $\Sr\big(t^{-1}\gll[t^{-1}]\big)$. Furthermore, we will only consider
the elements $\ol\Theta_m^{\tss(r)}$; the extension to the Pfaffian-type invariants
in type $D_n$ will then be obvious.

Using matrices \eqref{matr},
introduce power series in a variable $z$ by
\beql{frz}
\Fc^{(r)}(z)=\sum_{p=0}^{\infty} F^{(r)}[-p-1]\tss z^p
\eeq
and write the series
\beql{thes}
\ol\Theta_m(z,u)=\tr\ts\big(\Fc^{(0)}(z)+\dots+\Fc^{(\ell)}(z)\tss u^{\ell}\big)^m
\eeq
as a polynomial in $u$,
\ben
\ol\Theta_m(z,u)=\sum_r \ol\Theta^{\ts(r)}_m(z)\tss u^r.
\een
The element $T^s\tss\ol\Theta_m^{\ts(r)}$ then equals $s!$ times
the coefficient of $z^s$ in the series $\ol\Theta^{\ts(r)}_m(z)$.

On the other hand, we can use the
non-degenerate invariant
symmetric bilinear form on the Lie algebra $\gll$ defined by
\ben
\big\langle X_0+X_1v+\dots+X_{\ell}v^{\ell},Y_0+Y_1v+\dots+Y_{\ell}v^{\ell}\big\rangle^{}_{\ell}=
\langle X_0,Y_{\ell}\rangle+\dots+\langle X_{\ell},Y_0\rangle
\een
to identify the $\gll[t]$-module
$t^{-1}\gll[t^{-1}]$ with the restricted dual $\gll[t]^*$.
Namely, we will use the isomorphism
such that
for any $X,Y\in \gll$ and all values $p,s\geqslant 0$
the image of the
element $X[-p-1]\in t^{-1}\gll[t^{-1}]$ is the functional
whose value on $Y[s]\in \gll[t]$ equals $\de_{p\tss s}\langle X,Y\rangle^{}_{\ell}$.
Denote the image of the matrix $F^{(r)}[-p-1]$ under this isomorphism
by $\Xc^{(\ell-r)}[p]$, and the corresponding image of the series \eqref{frz}
by $\Xc^{(\ell-r)}(z)$. Accordingly, the image of the series \eqref{thes}
will then be written as
\ben
\Tc_m(z,u)=\tr\ts\big(\Xc^{(\ell)}(z)+\dots+\Xc^{(0)}(z)\tss u^{\ell}\big)^m
\een
so that
\ben
u^{-m\ell}\ts \Tc_m(z,u)=\tr\ts\big(\Xc^{(0)}(z)+\dots+\Xc^{(\ell)}(z)\tss u^{-\ell}\big)^m.
\een
Expand this as a polynomial in $u^{-1}$
\ben
\tr\ts\big(\Xc^{(0)}(z)+\dots+\Xc^{(\ell)}(z)\tss u^{-\ell}\big)^m
=\sum_{r\geqslant 0} \Tc_m^{(r)}(z)\tss u^{-r}
\een
thus defining the power series
\ben
\Tc_m^{(r)}(z)=\sum_{s=0}^{\infty}\Tc_m^{(r)}[s]\ts z^s.
\een

By identifying $\Sr(\gll[t]^*)$ with
the algebra of polynomial functions
$\Fun \gll[t]$ on $\gll[t]$ we come to proving that
all coefficients of the series $\Tc_m^{(r)}(z)$
with $m$ running over the values specified in Table~1
and $r=0,1,\dots,\ell$ are algebraically independent
generators of the algebra of $\gll[t]$-invariants of $\Fun \gll[t]$.
It is sufficient to establish the corresponding property for each truncated
algebra $\gllm$ defined as the quotient of $\gll[t]$ by the ideal
$t^{\ell'+1}\gll[t]$. More precisely, we need to
show that for any nonnegative integer $\ell'$, the coefficients
$\Tc_m^{(r)}[s]$ with the same conditions on $m$ and $r$ as above, and with
$s$ taking the values $0,1,\dots,\ell'$, are algebraically independent
generators of the algebra of $\gllm$-invariants of $\Fun \gllm$.
Writing the definition of the coefficients
$\Tc_m^{(r)}[s]$ in the form
\ben
\tr\ts\Big(\sum_{r=0}^{\ell}\sum_{p=0}^{\ell'}\Xc^{(r)}[p] z^pu^{-r}\Big)^m
=\sum_{r,s\geqslant 0}\Tc_m^{(r)}[s] z^su^{-r},
\een
observe that the required
property of these coefficients holds by
\cite[Theorem~5.4(b)]{ms:ip}. This follows since
the traces $\tr \ts \Xc^m$ with $m$ running over the values specified in Table~1,
are algebraically independent
generators of the algebra of $\g$-invariants of $\Fun \g$, where $\Xc$ denotes the image
of the matrix $F$ under the isomorphism $\Sr(\g)\cong \Fun \g$
defined by the form \eqref{killi}.
\epf

\section*{Acknowledgements}

I am grateful to Oksana Yakimova for providing
very useful answers to my questions.
The support of
the Australian Research Council, grant DP180101825
is acknowledged.

\section*{Data Availability}

The data that supports the findings of this study are available within the article.

\end{document}